\documentclass[a4paper,10pt]{amsart}
\usepackage{amssymb}
\numberwithin{equation}{section}
\newtheorem{theorem}{Theorem}[section]

\newtheorem{lem}[theorem]{Lemma}

\theoremstyle{remark}
\newtheorem{rem}[theorem]{Remark}

\renewcommand{\hat}{\widehat}

\author[J.~Benameur]{Jamel Benameur}
\address{Department of Mathematics, College of Science, King Saud University\\
Riyadh 11451, Kingdom of Saudi Arabia}
\email{\sl jbenameur@ksu.edu.sa}
\author[L.~Jlali]{Lotfi Jlali}
\address{Department of Mathematics, College of Science, King Saud University\\
Riyadh 11451, Kingdom of Saudi Arabia}
\email{\sl ljlali@ksu.edu.sa}
\title[Long time decay for 3D-NSE in Gevrey-Sobolev spaces]
{Long time decay for 3D-NSE  in Gevrey-Sobolev spaces}
\date{\today}
\begin{document}
\begin{abstract}In this paper we  prove, if $u$ is a global solution to Navier-Stokes equations in the Sobolev-Gevrey spaces $H^1_{a,\sigma}(\mathbb R^3)$, then $\|u(t)\|_{H^1_{a,\sigma}}$ decays to zero as time goes to infinity. Fourier analysis is used.
\end{abstract}

%@@@@@@@@@@@@@@@@@@@@@@@@@@@@@@@@@@@@%@@@@@@@@@@@@@@@@@@@@@@@@@@@@@@@@@@@@%@@@@@@@@@@@@@@

\subjclass[2000]{35-xx, 35Bxx, 35Lxx}
\keywords{Navier-Stokes Equations; Critical spaces; Long time decay}

%@@@@@@@@@@@@@@@@@@@@@@@@@@@@@@@@@@@@%@@@@@@@@@@@@@@@@@@@@@@@@@@@@@@@@@@@@%@@@@@@@@@@@@@@
\maketitle
\tableofcontents

%@@@@@@@@@@@@@@@@@@@@@@@@@@@@@@@@@@@@%@@@@@@@@@@@@@@@@@@@@@@@@@@@@@@@@@@@@%@@@@@@@@@@@@@@@

\section{Introduction}
The $3D$ incompressible Naviers-Stokes equations are given by:

$$
\left\{
  \begin{array}{lll}
     \partial_t u
 -\Delta u+ u.\nabla u  & =&\;\;-\nabla p\hbox{ in } \mathbb R^+\times \mathbb R^3\\
    {\rm div}\, u &=& 0 \hbox{ in } \mathbb R^+\times \mathbb R^3\\
    u(0,x) &=&u^0(x) \;\;\hbox{ in }\mathbb R^3,
  \end{array}
\right.
\leqno{(NSE)}
$$
where, we suppose that the fluid viscosity $\nu=1$, and  $u=u(t,x)=(u_1,u_2,u_3)$ and $p=p(t,x)$ denote respectively the unknown velocity and the unknown pressure of the fluid at the point $(t,x)\in \mathbb R^+\times \mathbb R^3$, $(u.\nabla u):=u_1\partial_1 u+u_2\partial_2 u+u_3\partial_3u$, and  $u^0=(u_1^o(x),u_2^o(x),u_3^o(x))$ is a given initial velocity. If $u^0$ is quite regular, the divergence free condition determines the pressure $p$.\\

We define the Sobolev-Gevrey spaces as follows; for $a, s\geq0$ and $\sigma>1$,
$$H^s_{a,\sigma}(\mathbb R^3)=\{f\in L^2(\mathbb R^3);\;e^{a|D|^{1/\sigma}}f\in H^s(\mathbb R^3)\}.$$
It is equipped with the norm
$$\|f\|_{H^s_{a,\sigma}}^2=\|e^{a|D|^{1/\sigma}}f\|_{H^s}$$
and its associated inner product
$$\langle f/g\rangle_{{H}^s_{a,\sigma}}=\langle e^{a|D|^{\frac{1}{\sigma}}}f/e^{a|D|^{\frac{1}{\sigma}}}g\rangle_{H^s}.$$

There are several authors who have studied the behavior of the norm of the solution to infinity in the different Banach spaces. For example:\\
Wiegner proved in \cite{MW1} that the $L^2$ norm of the solutions vanishes for any square integrable initial data, as times goes to infinity and gave a decay rate that seems to be optimal for a class of initial data. In \cite{MS2,MWS} M.E.Schonber and M.Wiegner derived some asymptotic properties of the solution and its higher derivatives under additional assumptions on the initial data. In \cite {BS1} J.Benameur and R.Selmi proved that if $u$ be a Leray solution of $(2d-NSE)$ then $\lim_{t\rightarrow\infty}\|u(t)\|_{L^2(\mathbb R^2)}=0$. In \cite{GIP} for the critical Sobolev spaces $\dot{H}^{\frac{1}{2}}$ I.Gallagher, D.Iftimie and F.Planchon
proved that $\|u(t)\|_{\dot{H}^{\frac{1}{2}}}$ goes to zero at infinity. In \cite{JB} J.Benameur proved if $u\in{\mathcal C}([0,\infty),{\mathcal{X}}^{-1}(\mathbb R^3))$ be a global solution to $3D$ Navier-Stokes equation, then $\|u(t)\|_{{\mathcal{X}}^{-1}}$ decay to zero as times goes to infinity.\\

We state our main result.
\begin{theorem}\label{theo1}
Let $a>0$ and $\sigma>1$. Let $u\in{\mathcal C}([0,\infty),H^1_{a,\sigma}(\mathbb R^3))$ be a global solution to $(NSE)$ system. Then
\begin{equation}\label{eq1}
\limsup_{t\rightarrow\infty}\|u(t)\|_{H^1_{a,\sigma}}=0.
\end{equation}
\end{theorem}
\begin{rem} The existence of local solutions to $(NSE)$ was studied in a recent paper \cite{BL1}.
\end{rem}
The paper is organized in the following way: In section $2$, we give some notations and important preliminary results. The section $3$ is devoted to prove that,  if $u\in{\mathcal C}(\mathbb R^+, H^1(\mathbb R^3))$ is a global solution to $(NSE)$, then $\|u(t)\|_{H^1}$ decays to zero as time goes to infinity. This proof uses the fact that
\begin{equation}\label{eq2}
\lim_{t\rightarrow\infty}\|u(t)\|_{\dot{H}^{\frac{1}{2}}}= 0
\end{equation}
and the energy estimate
\begin{eqnarray}\label{enq1}
\|u(t)\|_{L^2}^2+\int_0^t\|\nabla u(\tau)\|_{L^2}^2d\tau\leq\|u^0\|_{L^2}^2.
\end{eqnarray}
In section $4$, we generalize the results of Foias-Temam( see \cite{FT}) to $\mathbb R^3$ . In section $5$, we prove the main theorem. This proof is based on the obtained results in sections 3 and 4.

\section{Notations and preliminaries results}
\subsection{Notations}In this section, we collect some notations and definitions that will be used later.\\
$\bullet$ The Fourier transformation is normalized as
$$
\mathcal{F}(f)(\xi)=\widehat{f}(\xi)=\int_{\mathbb R^3}\exp(-ix.\xi)f(x)dx,\,\,\,\xi=(\xi_1,\xi_2,\xi_3)\in\mathbb R^3.
$$
$\bullet$ The inverse Fourier formula is
$$
\mathcal{F}^{-1}(g)(x)=(2\pi)^{-3}\int_{\mathbb R^3}\exp(i\xi.x)g(\xi)d\xi,\,\,\,x=(x_1,x_2,x_3)\in\mathbb R^3.
$$
$\bullet$ For $ s\in\mathbb R $, $H^s(\mathbb R^3)$ denotes the usual non-homogeneous Sobolev space on $\mathbb R^3$ and  $\langle./.\rangle_{H^s}$ denotes the usual scalar product on $H^s(\mathbb R^3)$.\\
$\bullet$ For $ s\in\mathbb R $, $\dot{H}^s(\mathbb R^3)$ denotes the usual homogeneous Sobolev space on $\mathbb R^3$ and $\langle./.\rangle_{\dot{H}^s}$ denotes the usual scalar product on $\dot{H}^s(\mathbb R^3)$.\\
$\bullet$ The convolution product of a suitable pair of functions $f$ and $g$ on $\mathbb R^3$ is given by
$$
(f\ast g)(x):=\int_{\mathbb R^3}f(y)g(x-y)dy.
$$
$\bullet$ If $f=(f_1,f_2,f_3)$ and $g=(g_1,g_2,g_3)$ are two vector fields, we set
$$
f\otimes g:=(g_1f,g_2f,g_3f),
$$
and
$$
{\rm div}\,(f\otimes g):=({\rm div}\,(g_1f),{\rm div}\,(g_2f),{\rm div}\,(g_3f)).
$$

\subsection{Preliminary results}
\begin{lem}\label{lem1}(See \cite{HB})
Let $(s,t)\in{\mathbb{R}^2}$, such that $s<{\frac{3}{2}}$ and $s+t>0$. Then, there exists a constant $C>0$, such that
$$
\|uv\|_{\dot{H}^{{s+t-{\frac{3}{2}}}}(\mathbb R^3)}\leq C(\|u\|_{\dot{H}^s(\mathbb R^3)}\|v\|_{\dot{H}^t(\mathbb R^3)}+\|u\|_{\dot{H}^t(\mathbb R^3)}\|v\|_{\dot{H}^s(\mathbb R^3)}).
$$
If $s<{\frac{3}{2}}$, $t<{\frac{3}{2}}$ and $s+t>0$,  then there exists a constant $C>0$, such that
$$
\|uv\|_{\dot{H}^{{s+t-{\frac{3}{2}}}}(\mathbb R^3)}\leq C \|u\|_{\dot{H}^s(\mathbb R^3)}\|v\|_{\dot{H}^t(\mathbb R^3)}.
$$
\end{lem}
\begin{lem}\label{lem3}
Let $f\in \dot H^{s_1}(\mathbb{R}^3)\cap \dot H^{s_2}(\mathbb{R}^3)$, where $s_1 < \frac{3}{2} < s_2 $. Then, there is a constant $c=c(s_1,s_2)$ such that
$$
\|f\|_{L^{\infty}(\mathbb{R}^3)}\leq \|\hat{f}\|_{L^{1}(\mathbb{R}^3)}\leq c\|f\|^{\frac{s_2-\frac{3}{2}}{s_2-s_1}}_{ \dot H^{s_1}(\mathbb{R}^3)}\|f\|^{\frac{\frac{3}{2}-s_1}{s_2-s_1}}_{\dot H^{s_2}(\mathbb{R}^3)}.
$$
\end{lem}
\noindent{\it Proof.} We have
\begin{eqnarray*}
\|f\|_{L^{\infty}}&\leq& \|\widehat{f}\|_{L^{1}}\\&\leq&  \int_{\xi} |\widehat{f(\xi)}|d\xi\\&\leq& \int_{|\xi|<{\lambda}}|\widehat{f(\xi)}|d\xi+\int_{|\xi|>{\lambda}}|\widehat{f(\xi)}|d\xi.
\end{eqnarray*}
We take
$$
I_1=\int_{|\xi|<{\lambda}}\frac{1}{|\xi|^{s_1}}|\xi|^{s_1}|\widehat{f(\xi)}|d\xi.
$$
Using the Cauchy-Schwarz inequality, we obtain
\begin{eqnarray*}
I_1&\leq&  \left(\int_{|\xi|<{\lambda}}\frac{1}{|\xi|^{s_1}}d\xi\right)^{\frac{1}{2}}\|f\|_ {\dot H^{s_1}(\mathbb{R}^3)}\\&\leq& \left(\int^{\lambda}_{0}\frac{1}{r^{2s_{1}-2}}dr\right)^{\frac{1}{2}}\|f\|_ {\dot H^{s_1}(\mathbb{R}^3)}\\&\leq& c_{s_1}\lambda^{\frac{3}{2}-s_1}\|f\|_ {\dot H^{s_1}(\mathbb{R}^3)}.
\end{eqnarray*}
Similarly, take
$$
I_2=\int_{|\xi|>{\lambda}}\frac{1}{|\xi|^{s_2}}|\xi|^{s_2}|\widehat{f(\xi)}|d\xi,
$$
we have
\begin{eqnarray*}
I_2&\leq&  \left(\int_{|\xi|>{\lambda}}\frac{1}{|\xi|^{s_2}}d\xi\right)^{\frac{1}{2}}\|f\|_ {\dot H^{s_2}}\\&\leq& \left(\int^{\infty}_{\lambda}\frac{1}{r^{2s_{2}-2}}dr\right)^{\frac{1}{2}}\|f\|_ {\dot H^{s_2}}\\&\leq&c_{s_2}\lambda^{\frac{3}{2}-s_2}\|f\|_{\dot H^{s_2}}.
\end{eqnarray*}
Therefore,
$$
\|f\|_{L^{\infty}}\leq A\lambda^{\frac{3}{2}-s_1}+B\lambda^{\frac{3}{2}-s_2}.
$$
with $A=c_{s_1}\|f\|_ {\dot H^{s_1}}$ and $B=c_{s_2}\|f\|_{\dot H^{s_2}}$.\\
Posing
$$
\varphi(\lambda)= A\lambda^{\frac{3}{2}-s_1}+B\lambda^{\frac{3}{2}-s_2}.
$$
Then, $\varphi'(\lambda)=0\Leftrightarrow\lambda=c(s_1,s_2)\left(\frac{B}{A}\right)^{\frac{1}{{s_2}-{s_1}}}$\\
So,
$$
||f||_{L^{\infty}(\mathbb{R}^3)}\leq c' A^{\frac{{s_2}-\frac{3}{2}}{{s_2}-{s_1}}} B^{\frac{{\frac{3}{2}}-{s_1}}{{s_2}-{s_1}}}.
$$
\hfill $\square$
\begin{rem}
In particular, for $s_1=1$ and $s_2=2$, where $f\in \dot H^{1}(\mathbb{R}^3)\cap \dot H^{2}(\mathbb{R}^3)$, we get
$$
||f||_{L^{\infty}}\leq ||f||^{\frac{1}{2}}_ {\dot H^{1}}||f||^{\frac{1}{2}}_ {\dot H^{2}}.
$$
\end{rem}

\section{Long time decay of $(NSE)$ system in $H^1(\mathbb R^3)$}
In this section, we want to prove: If $u\in{\mathcal C}(\mathbb R^+,H^1(\mathbb R^3))$ is a global solution to $(NSE)$ system, then
\begin{equation}\label{eq3}
\limsup_{t\rightarrow\infty}\|u(t)\|_{H^1}=0.
\end{equation}
This proof is done in two steps.\\\\
$\bullet$ Step 1: In this step, we shall prove that
\begin{equation}\label{eq4}
\limsup_{t\rightarrow\infty}\|u(t)\|_{\dot H^1}=0.
\end{equation}
We have
$$
\frac{1}{2}\frac{d}{dt}\|u\|_{\dot H^{\frac{1}{2}}}^2+ \|u\|_{\dot H^{\frac{3}{2}}}^2\leq c \|u\|_{\dot{H}^{\frac{1}{2}}} \|u\|_{\dot {H}^{\frac{3}{2}}}^2.
$$
From (\ref{eq2}), let $t_0>0$ such that $\|u(t_0)\|_{\dot H^{\frac{1}{2}}}<\frac{1}{2c}$. Then
$$
\frac{1}{2}\frac{d}{dt} \|u\|_{\dot H^{\frac{1}{2}}}^2+\frac{1}{2} \|u\|_{\dot H^{\frac{3}{2}}}^2\leq 0,\,\,\,\forall t\geq t_0.
$$
Integrating with respect to time, we obtain
\begin{eqnarray*}
\|u(t)\|_{\dot H^{\frac{1}{2}}}^2+\int_{t_0}^t \| u(\tau)\|_{\dot H^{\frac{3}{2}}}^2\leq  \|u(t_0)\|_{\dot H^{\frac{1}{2}}}^2,\,\,\,\,\forall t\geq t_0.
\end{eqnarray*}
Let $s>0$ and $c=c_s$. There exists $T_0=T_0(s,\nu,u^0)>0$, such that
$$
\|u(T_0)\|_{\dot H^{\frac{1}{2}}}<\frac{1}{2c_s}.
$$
Then
$$
\|u(t)\|_{\dot H^{\frac{1}{2}}}<c_s,\,\,\,\forall t\geq t_0.
$$
Now, for $s>0$ we have
\begin{eqnarray*}
\frac{1}{2}\frac{d}{dt}  \|u\|_{\dot {H}^s}^2+ \|u\|_{\dot H^{s+1}}^2&\leq& \|u\otimes u\|_{\dot {H}^s} \|u\|_{\dot H^{s+1}}\\&\leq&
c_s \|u\|_{\dot H^{\frac{1}{2}}}\|u\|_{\dot H^{s+1}}^2.
\end{eqnarray*}
Then
$$
\frac{1}{2}\frac{d}{dt}  \|u\|_{\dot {H}^s}^2+ \|u\|_{\dot H^{s+1}}^2\leq \frac{\nu}{2} \|u\|_{\dot H^{s+1}}^2,\,\,\,\,\forall t\geq T_0.
$$
Thus
$$
\frac{1}{2}\frac{d}{dt} \|u\|_{\dot {H}^s}^2+\frac{1}{2} \|u(t)\|_{\dot H^{s+1}}^2\leq 0,\,\,\, \forall t\geq T_0.
$$
So, for $T_0\leq t'\leq t$,
$$
\|u(t)\|_{\dot {H}^s}^2+\int_{t'}^t \|u(\tau)\|_{\dot H^{s+1}}^2d\tau \leq \|u(t')\|_{\dot {H}^s}^2.
$$
In particular, for $s=1$
\begin{eqnarray*}
\|u(t)\|_{\dot {H}^1}^2+\int_{t'}^t \|u(\tau)\|_{\dot H^{2}}^2 d\tau\leq \|u(t')\|_{\dot {H}^1}^2.
\end{eqnarray*}
Then $(t\rightarrow\|u(t)\|_{\dot H^1})$ is decreasing on $[T_0,\infty)$ and $u\in L^2([0,\infty),\dot H^2(\mathbb R^3))$.\\
Now, let $\varepsilon>0$ small enough. The $L^2$-energy estimate
$$
\|u(t)\|_{L^2}^2+2\int_{T_0}^t\|\nabla u(\tau)\|_{L^2}^2d\tau\leq\|u(T_0)\|_{L^2}^2,\;\;\forall t\geq T_0
$$
implies that $u\in L^2([T_0,\infty),\dot H^1(\mathbb R^3))$ and there is a time $t_{\varepsilon}\geq T_0$ such that
$$\|u(t_\varepsilon)\|_{\dot H^1}<\varepsilon.$$
As $(t\longrightarrow\|u(t)\|_{\dot H^1})$ is decreasing on $[T_0,\infty)$, then
$$\|u(t)\|_{\dot H^1}<\varepsilon,\,\,\forall t\geq t_\varepsilon.$$
Therefore (\ref{eq4}) is proved.
\hfill $\square$

\noindent$\bullet$ Step 2: In this step, we prove that
\begin{equation}\label{eq5}
\limsup_{t\rightarrow\infty}\|u(t)\|_{L^2}= 0.
\end{equation}
This proof is inspired by \cite{BB1} and \cite{BS1}. For $\delta>0$ and a given distribution $f$, we define the operators $A_{\delta}(D)$ and $B_{\delta}(D)$, as following:
$$
A_{\delta}(D)f=\mathcal{F}^{-1}(\mathbf{1}_{\{|\xi|<\delta\}}\mathcal{F}(f)),\,\,B_{\delta}(D)f=\mathcal{F}^{-1}(\mathbf{1}_{\{|\xi|\geq\delta\}
}\mathcal{F}(f)).
$$
It is clear that when applying $A_{\delta}(D)$ (respectively, $B_{\delta}(D)$) to any distribution, we are dealing with its low-frequency part(respectively, high- frequency part).\\
Let $u$ be a solution to $(NSE)$. Denote by $\omega_{\delta}$ and $\upsilon_{\delta}$, respectively, the low-frequency part and the high-frequency part of $u$ and so on ${\omega_{\delta}}^0$ and ${\upsilon_{\delta}}^0$ for the initial data $u^0$. Applying the pseudo-differential operators $A_{\delta}(D)$ to the $(NSE)$, we get
$$
\partial_t \omega_{\delta}-\nu\triangle\omega_{\delta}+A_{\delta}(D)\mathbb{P}(u.\nabla u)=0.
$$
Taking the $L^2(\mathbb R^3)$ inner product, we obtain
\begin{eqnarray*}
\frac{1}{2}\frac{d}{dt} \|\omega_{\delta}(t)\|_{L^2}^2+\|\nabla \omega_{\delta}(t)\|_{L^2}^2
&\leq&|\langle A_{\delta}(D)\mathbb{P}(u.\nabla u) /\omega_{\delta}(t) \rangle_{L^2}|\\&\leq& |\langle A_{\delta}(D)({\rm div}\,(u\otimes u)/\omega_{\delta}(t) \rangle_{L^2}|\\
&\leq&|\langle A_{\delta}(D)(u\otimes u ))/\nabla\omega_{\delta}(t) \rangle_{L^2}|\\&\leq& |\langle u\otimes u/\nabla\omega_{\delta}(t) \rangle_{L^2}|\\
&\leq&\|u\otimes u\|_{L^2}\|\nabla\omega_{\delta}(t)\|_{L^2}\\
&\leq&\|u\otimes u\|_{L^2}\|\nabla\omega_{\delta}(t)\|_{L^2}.
\end{eqnarray*}
Lemma \ref{lem1} yields
\begin{eqnarray*}
\frac{1}{2}\frac{d}{dt} \|\omega_{\delta}(t)\|_{L^2}^2+\|\nabla \omega_{\delta}(t)\|_{L^2}^2&\leq& C\|u(t)\| _{\dot H^{\frac{1}{2}}}
\|\nabla u(t)\|_{L^2}\|\nabla\omega_{\delta}(t)\|_{L^2}\\
&\leq& CM\|\nabla u(t)\|_{L^2}\|\nabla\omega_{\delta}(t)\|_{L^2}\;\;\;(M=\sup_{t\geq0}\|u(t)\| _{\dot H^{\frac{1}{2}}}).
\end{eqnarray*}
Integrating with respect to time, we obtain
$$
\|\omega_{\delta}(t)\|_{L^2}^2\leq \|{\omega_{\delta}}^0\|_{L^2}^2+CM\int_0^t \|\nabla u(\tau)\|_{L^2}\|\nabla\omega_{\delta}(\tau)\|_{L^2}d\tau.
$$
Hence, we have $\|\omega_{\delta}(t)\|_{L^2}^2\leq M_{\delta}$ for all $t\geq0$, where
$$
M_{\delta}=\|{\omega_{\delta}}^0\|_{L^2}^2+CM\int_0^{\infty} \|\nabla u(\tau)\|_{L^2}\|\nabla\omega_{\delta}(\tau)\|_{L^2}d\tau.
$$
On the one hand, it is clear that $\lim_{\delta\rightarrow0}\|{\omega_{\delta}}^0\|_{L^2(\mathbb R^3)}^2=0$.
On the other hand, the Lebesgue-Dominated Convergence Theorem implies that
\begin{equation}\label{eq6}
\lim_{\delta\rightarrow0}\int_0^{\infty}\|\nabla u(\tau)\|_{L^2}\|\nabla \omega_{\delta}(\tau)\|_{L^2}d\tau=0.
\end{equation}
Hence $\lim_{\delta\rightarrow0} M_{\delta}=0$, and thus
\begin{equation}\label{eq7}
\lim_{\delta\rightarrow0}\sup_{t\geq0}\|\omega_{\delta}(t)\|_{L^2}= 0.
\end{equation}
At this point, we note that it makes sense to take time equal to $\infty$ in the integral (\ref{eq6}). In fact, by definition of  $\omega_{\delta}$ we have
$\|\nabla \omega_{\delta}\|_{L^2}\leq\|\nabla u\|_{L^2}$ .\\
It is clear that, $\lim_{\delta\rightarrow0}\|\nabla \omega_{\delta}(t)\|_{L^2}=0$ almost everywhere. So, the integrand sequence
$$
\|\nabla u(t)\|_{L^2}\|\nabla \omega_\delta(t)\|_{L^2}
$$
converges point-wise to zero. Moreover, using the above computations and (\ref{enq1}), we obtain
$$
\|\nabla u(t)\|_{L^2}\|\nabla \omega_{\delta}(t)\|_{L^2}\leq \|\nabla u(t)\|_{L^2}^2\in L^1(\mathbb R^+).
$$
Thus, the integral sequence is dominated by an integrable function. Then the limiting function is integrable and one can take the time $T=\infty $ in (\ref{eq6}).\\
Now, let us investigate the high-frequency part. To do so, one applies the pseudo-differential operators $B_{\delta}(D)$ to the $(NSE)$ to get
$$
\partial_t \upsilon_{\delta}-\Delta\upsilon_{\delta}+B_{\delta}(D)\mathbb{P}(u.\nabla u)=0.
$$
Taking the Fourier transform with respect to the space variable, we obtain
\begin{eqnarray*}
\partial_t |\widehat{\upsilon_\delta}(t,\xi)|^2+2 |\xi|^2 |\widehat{\upsilon_\delta}(t,\xi)|^2&\leq&2|\mathcal{F}(B_{\delta}(D)\mathbb{P}(u.\nabla u)
)(t,\xi)| |\widehat{\upsilon_\delta}(t,\xi)|\\&\leq&2|\xi||\mathcal{F}(B_{\delta}(D)\mathbb{P}(u\otimes u))(t,\xi)|.|\widehat{\upsilon_\delta}(t,\xi)|\\&\leq&2|\mathcal{F}(u\otimes u)(t,\xi)|.|\widehat{\nabla\upsilon_\delta}(t,\xi)|.
\end{eqnarray*}
Multiplying the obtained equation by $\exp(2\nu|\xi|^2)$ and integrating with respect to time, we get
$$
|\widehat{\upsilon_\delta}(t,\xi)|^2\leq e^{-2 t|\xi|^2}|\widehat{\upsilon_\delta^0}(\xi)|^2+2\int_0^t e^{-2(t-\tau)|\xi|^2} |\mathcal{F}(u\otimes u)(\tau,\xi)|.|\widehat{\nabla\upsilon_\delta}(\tau,\xi)|d\tau.
$$
Since $|\xi|>\delta$, we have
$$
|\widehat{\upsilon_\delta}(t,\xi)|^2\leq e^{-2 t\delta^2}|\widehat{\upsilon_\delta^0}(\xi)|^2+2\int_0^t e^{-2(t-\tau)\delta^2} |\mathcal{F}(u\otimes u)(\tau,\xi)|.|\widehat{\nabla\upsilon_\delta}(\tau,\xi)|d\tau.
$$
Integrating with respect to the frequency variable $\xi$ and using Cauchy-Schwartz inequality, we obtain
$$
\|\upsilon_{\delta}(t)\|_{L^2}^2\leq  e^{-2 t{\delta}^2}\|\upsilon_{{\delta}^0}\|_{L^2}^2+2\int_0^t e^{-2(t-\tau){\delta}^2}
\|u\otimes u\|_{L^2}\|\nabla\upsilon_{\delta} \|_{L^2}d\tau.
$$
By the definition of $\upsilon_\delta$, we have
$$
\|\upsilon_{\delta}(t)\|_{L^2}^2\leq  e^{-2 t{\delta}^2}\|u^0\|_{L^2}^2+2\int_0^t e^{-2(t-\tau){\delta}^2}\|u\otimes u\|
_{L^2}\|\nabla u \|_{L^2}d\tau.
$$
Lemma \ref{lem1} and inequality (\ref{eq2}) yield
$$
\|\upsilon_{\delta}(t)\|_{L^2(\mathbb R^3)}^2\leq  e^{-2 t{\delta}^2}\|u^0\|_{L^2(\mathbb R^3)}^2+c\int_0^t e^{-2(t-\tau){\delta}^2}
\|u \|_{\dot H^{\frac{1}{2}}}\|\nabla u \|_{L^2}^2d\tau.
$$
$$\leq  e^{-2t{\delta}^2}\|u^0\|_{L^2}^2+CM\int_0^t e^{-2(t-\tau){\delta}^2}
\|\nabla u \|_{L^2}^2d\tau,\;\;(M=\sup_{t\geq0}\|u \|_{\dot H^{\frac{1}{2}}}).$$
Hence, $\|\upsilon_{\delta}(t)\|_{L^2}^2\leq N_{\delta}(t)$, where
$$
N_{\delta}(t)= e^{-2 t{\delta}^2}\|u^0\|_{L^2}^2+CM\int_0^{\infty} e^{-2(t-\tau){\delta}^2}\|\nabla u \|_{L^2}^2d\tau.
$$
Using Young inequality and inequality (\ref{enq1}), we get $N_{\delta}\in L^1(\mathbb R^+)$ and
$$
\int_0^{\infty}N_{\delta}(t)dt\leq \frac{\|u^0\|_{L^2}^2}{2\delta^2}+\frac{CM\|u^0\|_{L^2}^2}{4\delta^2}.
$$
So $t\rightarrow\|\upsilon_{\delta}(t)\|_{L^2}^2$ is continuous and belongs to $L^1(\mathbb R^+)$.\\
Now, let $\varepsilon>0$. At first, ( \ref{eq7}) implies that there exist some $\delta_0>0$ such that
\begin{eqnarray*}
\|\omega_{\delta_0}(t)\|_{L^2}\leq \varepsilon/2,\,\forall\,t\geq0.
\end{eqnarray*}
Let us consider the set $\mathrm{R}_{\delta_0}$ defined by $\mathrm{R}_{\delta_0}:=\{t\geq0,\,\|\upsilon_{\delta}(t)\|_{L^2(\mathbb R^3)}>\varepsilon/2\}$. If we denote by $\lambda_1(\mathrm{R}_{\delta_0})$ the Lebesgue measure of $\mathrm{R}_{\delta_0}$, we have
$$
\int_0^{\infty}\|\upsilon_{\delta_0}(t)\|_{L^2(\mathbb R^3)}^2dt\geq\int_{\mathrm{R}_{\delta_0}}\|\upsilon_{\delta}(t)\|_{L^2(\mathbb R^3)}^2dt\geq (\varepsilon/2)^2 \lambda_1(\mathrm{R}_{\delta_0}).
$$
By doing this, we can deduce that $\lambda_1(\mathrm{R}_{\delta_0})= T^\varepsilon_{\delta^0}<\infty$, and there exists $t^\varepsilon_{\delta^0}>T^\varepsilon_{\delta^0}$ such that
\begin{eqnarray*}
\|\upsilon_{\delta_0}(t^\varepsilon_{\delta^0})\|_{L^2}^2\leq (\varepsilon/2)^2.
\end{eqnarray*}
So, $\|u(t^\varepsilon_{\delta^0})\|_{L^2}\leq \varepsilon$ and from  (\ref{enq1}) we have
$$
\|u(t)\|_{L^2}\leq \varepsilon,\,\,\,\forall t\geq t^\varepsilon_{\delta^0}.
$$
This completes the proof of (\ref{eq5}). \hfill $\square$

\section{Generalization of  Foias-Temam result in $H^1(\mathbb R^3)$}
In \cite{FT} Fioas and Teamam  proved an analytic property for the Navier-Stokes equations on the torus $\mathbb T^3=\mathbb R^3/\mathbb Z^3$. Here, we give a similar result on whole space $\mathbb R^3$.
\begin{theorem}\label{theo2}
We assume that $u^0\in H^{1}(\mathbb{R}^3)$.  Then, there exists a time $T$ that depends only on the $\|u^0\|_{\dot H^{1}(\mathbb{R}^3)}$, such that:\\
$(NSE)$ possesses on $(0,T)$ a unique regular solution $u$ such that $(t\rightarrow e^{\nu t|D|}u(t))$ is continuous from $[0,T] $ into $H^{1}(\mathbb{R}^3)$). Moreover if $u\in{\mathcal C}(\mathbb R^+,H^1(\mathbb R^3))$ is a global solution to $(NSE)$ system, then there are $M\geq0$ and $t_0>0$ such that
$$\|e^{t_0|D|}u(t)\|_{H^1(\mathbb{R}^3)}\leq M,\;\;\forall t\geq t_0.$$
\end{theorem}

Before proving this theorem, we need the following lemmas
\begin{lem}\label{lem5}
Let $ t\mapsto e^{t|D|}u\in H^{2}(\mathbb{R}^3)$, where $|D|=(\Delta)^{\frac{1}{2}} $. Then
$$
\|e^{t|D|}u.\nabla v\|_{L^2(\mathbb{R}^3)}\leq \|e^{t|D|}u\|^{\frac{1}{2}}_{H^1(\mathbb{R}^3)} \|e^{t|D|} u\|^{\frac{1}{2}}_{H^2(\mathbb{R}^3)} \|e^{t|D|}\Delta^{\frac{1}{2}}v\|_{L^2(\mathbb{R}^3)}.
$$
\end{lem}
\noindent{\it Proof.} We have
\begin{eqnarray*}
\|e^{t|D|}u.\nabla v\|_{L^2}&=&\int_{\mathbb{R}^3} e^{2t|\xi|}|\widehat{u.\nabla v}(\xi)|^2d\xi\\
&\leq& \int_{\mathbb{R}^3} e^{2t|\xi|}\left(\int_{\mathbb{R}^3} |\hat{u}(\xi-\eta)||\widehat{\nabla v}(\eta|d\eta\right)^2d\xi\\
&\leq& \int_{\mathbb{R}^3}\left(\int_{\mathbb{R}^3} e^{t|\xi|}|\hat{u}(\xi-\eta)||\widehat{\nabla v}(\eta|d\eta\right)^2d\xi\\
&\leq&\int_{\mathbb{R}^3}\left(\int_{\mathbb{R}^3}\left(e^{t|\xi-\eta|}|\hat{u}(\xi-\eta)|\right)
\left(e^{t|\eta|}|\eta||\hat{v}(\eta)|\right)d\eta\right)^2d\xi\\
&\leq&\left(\int_{\mathbb{R}^3} e^{t|\xi|}|\hat{u}(\xi)|d\xi\right)^2||e^{t|D|}\Delta^{\frac{1}{2}}v||_{L^2}.
\end{eqnarray*}
Hence, for $f={\mathcal F}^{-1}( e^{t|\xi|}|\hat{u}(\xi)|)\in H^{2}(\mathbb{R}^3)$ and $(s_1=1;\;s_2=2)$, lemma \ref{lem3} gives the desired result.\\\\
\begin{lem}\label{lem6}
Let $t\mapsto e^{t|D|}u \in H^{2}(\mathbb{R}^3)$. Then
$$
\left|\langle e^{t|D|}(u.\nabla v) /e^{t|D|}w \rangle_{{H^1}}\right|\leq\|e^{t|D|}u\|^{\frac{1}{2}}_{H^1} \|e^{t|D|} u\|^{\frac{1}{2}}_{H^2} \|e^{t|D|}\Delta^{\frac{1}{2}}v\|_{L^2} \|e^{t|D|}\Delta w\|_{L^2}.
$$
\end{lem}
\noindent{\it Proof.}\\
We have
\begin{eqnarray*}
\langle u.\nabla v / w \rangle_{H^1}=\sum_{|j|=1}\langle\partial_j ( u.\nabla v)/\partial_j w\rangle_{L^2}\\=- \sum_{|j|=1}\langle u.\nabla v/\partial^{2}_j w\rangle_{L^2}\\=-\sum_{|j|=1}\langle u.\nabla v/ \Delta w\rangle_{L^2}.
\end{eqnarray*}
Then
\begin{eqnarray*}
\left|\langle e^{t|D|}u.\nabla v /e^{t|D|}w \rangle_{H^1}\right|=\left|\langle e^{t|D|}u.\nabla v /e^{t|D|}\Delta w \rangle_{L^2}\right|\\\leq\|e^{t|D|}u.\nabla v\|_{L^2} \|e^{t|D|}\Delta w\|_{L^2}
\end{eqnarray*}
Finally, using lemma \ref{lem5}, we obtain the desired result.\\\\
\noindent{\it Proof of theorem \ref{theo2}.} We have
\begin{equation}\label{eq8}
\partial_{t}u -\Delta u +u.\nabla u =-\nabla p.
\end{equation}
Applying  the fourier transform to the last equation  and multiplying by $\overline{\widehat{u}}$, we have
\begin{eqnarray}\label{eq9}
\partial_{t} \widehat{u}.\overline{\widehat{u}} + |\xi|^2 |\widehat{u}|^2 = -(\widehat{u.\nabla u}). \overline{\widehat{u}}.
\end{eqnarray}
Again, the fourier (bar) of (\ref{eq8})  multiplied by $\widehat{u}$ gives
\begin{eqnarray}\label{eq10}
\partial_{t} \overline{\widehat{u}}.\widehat{u}+  |\xi|^2 |\widehat{u}|^2 = -( \overline{\widehat{u.\nabla u}}). \widehat{u}.
\end{eqnarray}
Hence, the sum of (\ref{eq9}) and  (\ref{eq10}) yields
\begin{eqnarray*}
\partial_{t}|\widehat{u}|^2+ 2 |\xi|^2 |\widehat{u}|^2 = - 2 Re ((\widehat{u.\nabla u}) . \widehat{u}).
\end{eqnarray*}
This implies
$$
\partial_{t}|\widehat{u}|^2 (1+|\xi|^2) e^{2t|\xi|}+ 2(1+|\xi|^2) |\xi|^2 e^{2t|\xi|} |\widehat{u}|^2 = - 2 Re ((\widehat{u.\nabla u}) . \widehat{u})(1+|\xi|^2) e^{2t|\xi|}.
$$
Then
$$
\int_{\mathbb{R}^3}(1+|\xi|^2) e^{2t|\xi|}\partial_{t}|\widehat{u}(\xi)|^2 d\xi +2 \int_{\mathbb{R}^3}(1+|\xi|^2)|\xi|^2 e^{2t|\xi|} |\widehat{u}(\xi)|^2 d\xi = - 2 Re \int_{\mathbb{R}^3}((\widehat{u.\nabla u}) . \widehat{u})(1+|\xi|^2) e^{2t|\xi|} d\xi.
$$
Thus
$$
\langle e^{t|D|}\partial_{t} u / e^{t|D|} u \rangle_{H^1} + 2 \|e^{t|D|}\nabla u\|^2_{H^1(\mathbb{R}^3)} = - 2 Re \langle e^{t|D|} (u.\nabla u) /e^{t|D|}u \rangle_{H^1}.
$$
At time $\tau$, we have
\begin{eqnarray}\label{enq2}
\langle e^{\tau|D|} u'(\tau)/e^{\tau|D|} u(\tau) \rangle_{H^1} + 2 \|e^{\tau|D|}\nabla u\|^2_{H^1} = - 2 Re \langle e^{\tau|D|} (u.\nabla u)/e^{t|D|}u \rangle_{H^1}.
\end{eqnarray}
Therefore
\begin{eqnarray*}
\langle e^{t|D|} u'(t)/ e^{t|D|} u(t) \rangle_{H^1}&=& \langle ( e^{t|D|} u(t))'-|D|e^{t|D|}u(t) / e^{t|D|} u(t) \rangle_{H^1}\\&=& \frac{1}{2} \frac{d}{dt} \|e^{t|D|}u\|^2_{H^1}-\langle e^{t|D|}|D| u(t) / e^{t|D|} u(t) \rangle_{H^1}\\&\geq& \frac{1}{2} \frac{d}{dt} \|e^{t|D|}u\|^2_{H^1}- \|e^{t|D|}u\|_{H^1} \|e^{t|D|}u\|_{H^2}.
\end{eqnarray*}
Using the Young inequality, we obtain
\begin{eqnarray}\label{enq3}
\frac{d}{dt} \|e^{t|D|}u\|^2_{H^1} -2 \|e^{t|D|}u\|^2_{H^1} -\frac{1}{2}\|e^{t|D|}u\|^2_{H^2}\leq 2 \langle e^{t|D|} u'(t) / e^{t|D|} u(t) \rangle_{H^1}.
\end{eqnarray}
Hence, using the lemma \ref{lem6} and Young inequality the right hand of (\ref{enq2}) satisfies
\begin{eqnarray*}
|-2 Re\langle e^{\tau|D|} u\nabla u / e^{\tau|D|} u \rangle_{H^1}|&\leq&2 \|e^{\tau|D|} u\|^{\frac{1}{2}}_{H^1} \|e^{\tau|D|} u\|^{\frac{1}{2}}_{H^1} \|e^{\tau|D|} |D| u\|_{L^2} \|e^{\tau|D|}\Delta u\|_{L^2}\\&\leq&2 \|e^{\tau|D|}u\|^{\frac{3}{2}}_{H^1} \|e^{\tau|D|}u\|^{\frac{3}{2}}_{H^2}\\&\leq& \frac{3}{4}\|e^{\tau|D|}u\|^2_{H^2} + \frac{c_1}{2} \|e^{\tau|D|}u\|^6_{H^1},
\end{eqnarray*}
where $c_1$ is a positive  constant.\\
Then  (\ref{enq2}) yields
\begin{eqnarray}\label{enq4}
\langle e^{t|D|}u'(t) / e^{t|D|} u(t) \rangle_{H^1}+2 \|e^{t |D|} \nabla u\|^2_{H^1}&\leq& \frac{3}{4}\|e^{t|D|}u\|^2_{H^2}+\frac{c_1}{2}\|e^{t|D|}u\|^6_{H^1}.
\end{eqnarray}
Hence, using  (\ref{enq3}) and (\ref{enq4}), we get
\begin{eqnarray*}
\frac{d}{dt} \|e^{t|D|}u\|^2_{H^1} + 2 \|e^{t |D|} \nabla u\|^2_{H^1}&\leq& 4 \|e^{t|D|}u\|^2_{H^1} + c_1\|e^{t|D|}u\|^6_{H^1}\\&\leq& c_2 + 2 c_1\|e^{t|D|}u\|^6_{H^1},
\end{eqnarray*}
where also $ c_2$ is a positive constant.\\
Finally, we obtain
\begin{eqnarray*}
y'(t)\leq K_1 y^3(t),
\end{eqnarray*}
where
\begin{equation*}
 y(t)=1 + \|e^{t|D|} u(t) \|^2_{H^1} \quad\mbox{and}\quad K_1 =  2 c_1+c_2.
\end{equation*}
Then
\begin{eqnarray*}
y(t)\leq y(0)+K_1 \int^t_0y^3(s)ds.
\end{eqnarray*}
Let
\begin{equation*}
T_1=\frac{2}{K_1 y^2(0)}
\end{equation*}
and $o < T \leq T^*$ such that\,\,$T =\sup \{t\in[0,T^*)\,|\,\sup_{0\leq s\leq t}y(s)\leq 2y(0)\}$. Hence for $0 \leq t \leq\min(T_1,T)$, we have
\begin{eqnarray*}
y(t)&\leq& y(0)+K_1 \int^t_0 y^3(s)ds\\&\leq& y(0)+K_1 \int^t_0 8 y^3(0)ds\\&\leq& \left(1+K_1 8 T_1 y^2(0)\right)y(0).
\end{eqnarray*}
Taking $1+K_1 8 T_1 y^2(0)<2$, we get $ T>T_1$. Then
$$
y(t)\leq2y(0),\,\,\forall t \in [0,T_1].
$$
Therefore $ t\mapsto e^{t|D|}u(t)\in H^{1(\mathbb{R}^3)},\,\,\forall t \in [0,T_1]$.\\
In particular
$$
\| e^{T_1|D|}u(T_1)\|^2_{H^1} \leq 2+2\|u_0\|^2_{H^1}.
$$
Now, if we know that
\begin{eqnarray*}
\|u(t)\|_{H^1}\leq M_1 \,\,\,   \forall  t \geq 0.
\end{eqnarray*}
Defining the system
$$
\left\{
  \begin{array}{lll}
     \partial_t w
 -\Delta w+ w.\nabla w  & =&\;\;-\nabla p_2\hbox{ in } \mathbb R^+\times \mathbb R^3,\\
    {\rm div}\, w &=& 0 \hbox{ in } \mathbb R^+\times \mathbb R^3,\\
    w(0) &=&u(b) \;\;\hbox{ in }\mathbb R^3,
  \end{array}
\right.
\leqno{}
$$
where $w(t)=u(T+t)$.\\
Using a similar technic, we can prove that there exists $T_2 =\frac{2}{K_1}(1+M^{2}_1)^{-2}$  such that
$$
y(t)=1+\|e^{t|D|}w(t)\|^2_{H^1} \leq 2(1+M^{2}_1), \,\,\, \forall t \in [0,T_2].
$$
This implies that $1+\|e^{t|D|}u(T+t)\|^2_{H^1} \leq 2(1+M^{2}_1)$. Hence, for $t=T_2$ we have
$$
\|e^{T_2|D|}u(T+T_2)\|^2_{H^1} \leq 2(1+M^{2}_1).
$$
Since $t=T+T_2\geq T_2, \,\,\ \forall T\geq 0$, we obtain
$$
\|e^{T_2|D|}u(t)\|^2_{H^1} \leq 2(1+M^{2}_1),\,\,\, \forall t\geq T_2.
$$
Then
$$
\|e^{T_2|D|}u(t)\|^2_{H^1} \leq 2(1+M^{2}_1),\,\,\, \forall t\geq T_2,
$$
where
\begin{equation*}
T_2 = T_2(M_1) =\frac{2}{K_1}(1+M^{2}_1)^{-2}.
\end{equation*}
\hfill $\square$\\

\section{Proof of main result}
In this section, we prove the main theorem \ref{theo1}. This proof uses the result of sections $3$ and $4$.\\
 Let $u\in{\mathcal C}(\mathbb R^+, H^1_{a,\sigma}(\mathbb R^3) )$.
As $ H^1_{a,\sigma}(\mathbb R^3)\hookrightarrow H^1(\mathbb R^3)$, then $u\in{\mathcal C}(\mathbb R^+,H^1(\mathbb R^3))$.\\
Applying the theorem \ref{theo2}, there exist $t_0>$ and $\alpha>0$ such that
\begin{eqnarray}\label{enq4}
\|e^{\alpha|D|}u(t)\|_{H^1}\leq c_0=2+M_1^2,\;\;\;\forall t\geq t_0,
\end{eqnarray}
where $\alpha=\varphi(t_0)$ and $t_0=\frac{2}{K_1}(1+M_1^2)^{-2}$.\\
Therefore, let $a>0$, $\beta>0$. It shows that there exists $c_3\geq0$ such that
\begin{eqnarray*}
ax^{\frac{1}{\sigma}}\leq c_3+\beta x,\,\,\,\forall x\geq0.
\end{eqnarray*}
Indeed; $\frac{1}{\sigma}+\frac{\sigma-1}{\sigma}=\frac{1}{p}+\frac{1}{q}=1$. Using the Young inequality, we obtain
\begin{eqnarray*}
ax^{\frac{1}{\sigma}}&=&a \beta^{\frac{-1}{\sigma}}(\beta^{\frac{1}{\sigma}}x^{\frac{1}{\sigma}})\\&\leq&\frac{(a \beta^{\frac{-1}{\sigma}})^q}{q}+\frac{( \beta^{\frac{1}{\sigma}}x^{\frac{1}{\sigma}})^p}{p}\\&\leq&c_3+\frac{\beta x}{\sigma}\\&\leq&c_3+\beta x,
\end{eqnarray*}
where $c_3=\frac{\sigma-1}{\sigma}a^{\frac{\sigma}{\sigma-1}}\beta^{\frac{1}{1-\sigma}}$.\\
Take $\beta=\frac{\alpha}{2}$, using (\ref{enq4}) and the  Cauchy Schwarz inequality, we have
\begin{eqnarray*}
\|u(t)\|_{ H^1_{a,\sigma}}&=&\|e^{a|D|^{1/\sigma}}u(t)\|_{H^1}\\&=& \int(1+|\xi|^2) e^{2a|\xi|^{1/\sigma}}|\widehat{u}(t,\xi)|^2d\xi\\&=&\int(1+|\xi|^2) e^{2(c_3+\beta|\xi|})|\widehat{u}(t,\xi)|^2d\xi\\&=&\int(1+|\xi|^2) e^{2c_3}e^{\alpha |\xi|}|\widehat{u}(t,\xi)|^2d\xi\\&\leq&e^{2c_3}\left(\int(1+|\xi|^2)|\widehat{u}(t,\xi)|^2d\xi\right)^{\frac{1}{2}}\left(\int(1+|\xi|^2)e^{2\alpha|\xi|}|\widehat{u}(t,\xi)|^2d\xi\right)^{\frac{1}{2}}\\&\leq&e^{2c_3}\|u\|_{ H^1}^{\frac{1}{2}}\|e^{a|D|}u(t)\|_{H^1}^{\frac{1}{2}}\\&\leq&c\|u\|_{ H^1}^{\frac{1}{2}},
\end{eqnarray*}
where $c=e^{2c_3}c_0^{\frac{1}{2}}$.\\
Using (\ref{eq3}), we get
$$
\limsup_{t\rightarrow\infty}\|e^{a|D|^{1/\sigma}}u(t)\|_{H^1}=0.
$$
\hfill $\square$

\end{document}